January 02, 2017

# Binomial transform of products

Khristo N. Boyadzhiev
Department of Mathematics and Statistics,
Ohio Northern University, Ada, OH 45810 USA
k-boyadzhiev@onu.edu

**Abstract**. Given the binomial transforms $\{b_n\}$ and $\{c_n\}$ of the sequences $\{a_n\}$ and $\{d_n\}$ correspondingly, we compute the binomial transform of the sequence $\{a_n c_n\}$ in terms of $\{b_n\}$ and $\{d_n\}$. In particular, we compute the binomial transform of the sequences $\{n(n-1)...(n-m+1)a_n\}$ and $\{a_k x^k\}$ in terms of $\{b_n\}$. Further applications include new binomial identities with the binomial transforms of the products $H_n B_n$, $H_n F_n$, $H_n L_n(x)$, and $B_n F_n$, where $H_n$, $B_n$, $F_n$, and $L_n(x)$ are correspondingly the harmonic numbers, the Bernoulli numbers, the Fibonacci numbers, and the Laguerre polynomials.

## 1. Introduction and main results

Given a sequence $\{a_n\}$, its binomial transform is the sequence $\{b_n\}$ defined by the formula

$$b_n = \sum_{k=0}^{n} \binom{n}{k} a_k \qquad (1.1)$$

with inversion

$$a_n = \sum_{k=0}^{n} \binom{n}{k} (-1)^{n-k} b_k \ .$$

The binomial transform is related to Euler's series transformation [2] and provides numerous nice and elegant binomial identities (see [1], [4], [7]). It is a powerful instrument in the theory of special numbers [6] and in combinatorics.

Our purpose is to develop a technique that helps to generate new binomial transform identities from old. When the binomial transform (1.1) is known, we want to compute the image sequence

$$\sum_{k=0}^{n}\binom{n}{k}a_k c_k$$

($n = 0,1,...$), where the transform

$$d_n = \sum_{k=0}^{n}\binom{n}{k}(-1)^{n-k}c_k, \quad c_n = \sum_{k=0}^{n}\binom{n}{k}d_k \qquad (1.2)$$

is also known. Work in this direction was started in the recent paper [1], where it was shown that

$$\sum_{k=1}^{n}\binom{n}{k}\frac{a_k}{k+\lambda} = n!\sum_{m=1}^{n}\frac{b_m}{m!(\lambda+m)(\lambda+m+1)...(\lambda+n)}$$

for any $\lambda \geq 0$. Writing for brevity $\nabla b_n = b_n - b_{n-1}, \nabla^0 b_n = b_n$ for the backward difference, we have also the following result from [1]

$$\sum_{k=0}^{n}\binom{n}{k}k^p a_k = (n\nabla)^p b_n \qquad (1.3)$$

for any $0 \leq p \leq n$. This will be needed later.

In this paper we prove the following theorem.

**Theorem 1**. Let $\{a_n\}$ and $\{c_n\}$ be two sequences and let $\{b_n\}$ and $\{d_n\}$ be defined by (1.1) and (1.2). Then we have the identity

$$\sum_{k=0}^{n}\binom{n}{k}a_k c_k = \sum_{m=0}^{n}\binom{n}{m}d_m \nabla^m b_n, \qquad (1.4)$$

where $\nabla b_n = b_n - b_{n-1}$ with $\nabla^0 b_n = b_n$.

For the symmetric version of the transform $\{c_n\} \leftrightarrow \{d_n\}$, namely,

$$d_n = \sum_{k=0}^{n}\binom{n}{k}(-1)^k c_k, \quad c_n = \sum_{k=0}^{n}\binom{n}{k}(-1)^k d_k,$$

equation (1.4) takes the form

$$\sum_{k=0}^{n}\binom{n}{k}a_k c_k = \sum_{m=0}^{n}\binom{n}{m}(-1)^m d_m \nabla^m b_n . \qquad (1.5)$$

As we shall see, with appropriate choices of the sequences $\{a_n\}$ and $\{c_n\}$ this formula produces interesting new identities involving binomial polynomials and special numbers. In several cases the iterated differences $\nabla^m b_n$ can be computed explicitly.

The proof of the theorem is based on the lemma:

**Lemma 1**. Suppose the sequences $\{a_n\}$ and $\{b_n\}$ are defined from (1.1) and let $\nabla b_n = b_n - b_{n-1}$, $\nabla^0 b_n = b_n$. Then for every two integers $0 \leq m \leq n$ we have

$$\sum_{k=0}^{n}\binom{n}{k}k(k-1)\ldots(k-m+1)a_k = m!\binom{n}{m}\nabla^m b_n \qquad (1.6)$$

or, in a shorter form,

$$\sum_{k=0}^{n}\binom{n}{k}\binom{k}{m}a_k = \binom{n}{m}\nabla^m b_n . \qquad (1.7)$$

The case $m=1$ is (1.3) with $p=1$. It was proved in [1].

The next lemma presents one possible way to compute the RHS in the above identity.

**Lemma 2**. For any integers $0 \leq m \leq n$,

$$\binom{n}{m}\nabla^m b_n = \sum_{j=0}^{n}\binom{n}{j}\binom{j}{n-m}(-1)^{n-j}b_j . \qquad (1.8)$$

This can also be put in the form

$$\nabla^m b_n = m!\sum_{j=n-m}^{n}\frac{(-1)^{n-j}b_j}{(n-j)!(j-n+m)!} .$$

The proofs are given in section 5.

**Example 1**. Let $a_k = 1$ for all $k$. Then

$$b_n = \sum_{k=0}^{n}\binom{n}{k} = 2^n ,$$

and Lemma 1 gives for any $n = 0, 1, 2, \ldots$, and any $0 \leq m \leq n$,

$$\sum_{k=0}^{n}\binom{n}{k}k(k-1)...(k-m+1)=m!\binom{n}{m}\nabla^{m}2^{n}=m!\binom{n}{m}2^{n-m}, \qquad (1.9)$$

as by a simple computation we find $\nabla^{m}2^{n}=2^{n-m}$.

In the next section we apply our theorem for the case when $c_k = x^k$, and in section 3, among other things, we consider the case $a_k = (-1)^{k-1}H_k$, where $H_k$ are the harmonic numbers. We shall compute the iterated differences $\nabla^{m}H_n$ and prove the identity

$$\sum_{k=0}^{n}\binom{n}{k}(-1)^{k-1}H_k c_k = (-1)^{n-1}H_n d_n + \sum_{m=0}^{n-1}\frac{(-1)^{m}d_m}{n-m}, \qquad (1.10)$$

for any sequence $\{c_k\}$ where $c_k$ and $d_k$ are related by (1.2). In examples 13, 14, and 15 we apply this formula to the cases where $c_k$ are correspondingly the Fibonacci numbers, the Bernoulli numbers, and the Lagurre polynomials.

A similar identity is proved also for the Fibonacci numbers,

$$\sum_{k=0}^{n}\binom{n}{k}(-1)^{k-1}F_k c_k = \sum_{m=0}^{n}\binom{n}{m}d_m F_{n-2m} \qquad (1.11)$$

with $c_k$ and $d_k$ as above.

## 2. Binomial polynomials

For a given sequence $\{a_n\}$ we consider the polynomials

$$p_n(x) = \sum_{k=0}^{n}\binom{n}{k}a_k x^k.$$

When the binomial transform (1.1) is known, we want to compute the polynomials $p_n(x)$ explicitly in terms of the numbers $p_n(1) = b_n$. We present here two solutions to this problem. They both follow from Theorem1 with different choices of the sequences $\{a_n\}$ and $\{c_n\}$.

**Corollary 1.** Assuming the binomial transform (1.1) is given, we have the identity

$$\sum_{k=0}^{n}\binom{n}{k}a_k x^k = \sum_{j=0}^{n}\binom{n}{j}b_j x^j(1-x)^{n-j} = (1-x)^n \sum_{j=0}^{n}\binom{n}{j}b_j\left(\frac{x}{1-x}\right)^{j}. \qquad (2.1)$$

*Proof.* In Theorem 1 we set

$$a_n = (-1)^n x^n = \sum_{k=0}^{n} \binom{n}{k} (-1)^{n-k} (1-x)^k$$

so that in view of (1.1), $b_n = (1-x)^n$. Simple computation shows that for any $z$, $m$, and $n$,

$$\nabla^m z^n = z^{n-m} (z-1)^m$$

and thus

$$\nabla^m b_n = \nabla^m (1-x)^n = (1-x)^{n-m} (-x)^m.$$

From Theorem 1

$$\sum_{k=0}^{n} \binom{n}{k} (-1)^k x^k c_k = \sum_{m=0}^{n} \binom{n}{m} d_m (-x)^m (1-x)^{n-m} = \sum_{m=0}^{n} \binom{n}{m} (-1)^m d_m x^m (1-x)^{n-m}.$$

Here we change notations in order to write this equation in terms of $\{a_n\}$ and $\{b_n\}$. Setting $a_k = (-1)^k c_k$ we have from (1.2)

$$(-1)^n d_n = \sum_{k=0}^{n} \binom{n}{k} (-1)^k c_k = \sum_{k=0}^{n} \binom{n}{k} a_k = b_n,$$

and the above equation becomes

$$\sum_{k=0}^{n} \binom{n}{k} a_k x^k = \sum_{m=0}^{n} \binom{n}{m} b_m x^m (1-x)^{n-m},$$

as needed.

*Second proof, independent of Theorem 1.* Using the inversion formula we can write

$$p_n(x) = \sum_{k=0}^{n} \binom{n}{k} a_k x^k = \sum_{k=0}^{n} \binom{n}{k} x^k \left\{ \sum_{j=0}^{k} \binom{k}{j} (-1)^{k-j} b_j \right\}$$

$$= \sum_{j=0}^{n} (-1)^j b_j \left\{ \sum_{k=j}^{n} \binom{n}{k} \binom{k}{j} (-1)^k x^k \right\},$$

and the rest follows from the well-known identity [4], (3.11.8)

$$\sum_{k=j}^{n}\binom{n}{k}\binom{k}{j}(-1)^k x^k = (-1)^j \binom{n}{j} x^j (1-x)^{n-j}. \tag{2.2}$$

**Remark 1**. Identity (2.2) itself follows from Corollary 1 applied to the convolution identity

$$\sum_{k=j}^{n}\binom{n}{k}\binom{k}{j}(-1)^k = (-1)^n \delta_{nj}.$$

Here are some examples of representations of the form (2.1).

**Example 2**. The generalized Stirling numbers $S(\alpha,n)$ of the second kind are defined by the binomial formula (see [3] and the references therein)

$$\sum_{k=0}^{n}\binom{n}{k}(-1)^k k^\alpha = (-1)^n n! S(\alpha,n), \tag{2.3}$$

where $\alpha$ is any complex number with $\mathrm{Re}\,\alpha > 0$. According to Corollary 1 we have

$$\sum_{k=0}^{n}\binom{n}{k}(-1)^k k^\alpha x^k = \sum_{j=0}^{n}\binom{n}{j} j! S(\alpha,j)(-1)^j x^j (1-x)^{n-j},$$

or, changing $x$ to $-x$,

$$\sum_{k=0}^{n}\binom{n}{k} k^\alpha x^k = \sum_{j=0}^{n}\binom{n}{j} j! S(\alpha,j) x^j (1+x)^{n-j} \tag{2.4}$$

When $\alpha$ is a non-negative integer, $S(\alpha, j)$ are the usual Stirling numbers of the second kind [6].

**Example 3**. Setting $x = 2$ in (2.1) we obtain the curious identity

$$\sum_{k=0}^{n}\binom{n}{k} 2^k a_k = \sum_{j=0}^{n}\binom{n}{j}(-1)^{n-j} 2^j b_j.$$

**Example 4**. With $x = \dfrac{1}{2}$ in (2.1),

$$\sum_{k=0}^{n}\binom{n}{k} 2^{n-k} a_k = \sum_{j=0}^{n}\binom{n}{j} b_j$$

which explains the action of the iterated binomial transform.

The representation (2.1) in the above corollary is short and simple, but its RHS is not a polynomial in standard form. From Theorem 1 we obtain also a second corollary:

**Corollary 2**. Suppose the sequence $\{b_n\}$ is the binomial transform of the sequence $\{a_n\}$. Then

$$\sum_{k=0}^{n}\binom{n}{k}a_k x^k = \sum_{m=0}^{n}\binom{n}{m}\nabla^m b_n (x-1)^m . \qquad (2.5)$$

*Proof.* Taking $c_k = x^k$ in Theorem 1, equation (1.4), we have the one-line proof

$$d_n = \sum_{k=0}^{n}\binom{n}{k}(-1)^{n-k} x^k = (-1)^n \sum_{k=0}^{n}\binom{n}{k}(-x)^k = (x-1)^n .$$

Formula (2.5) gives, in fact, the Taylor expansion of the polynomial $p_n(x)$ centered at $x=1$. Most existing examples of binomial transforms with "$x$" have this format.

**Example 5**. With $x=2$ in (2.5) we find

$$\sum_{k=0}^{n}\binom{n}{k}2^k a_k = \sum_{m=0}^{n}\binom{n}{m}\nabla^m b_n .$$

This together with Example 3 provide

$$\sum_{m=0}^{n}\binom{n}{m}\nabla^m b_n = \sum_{j=0}^{n}\binom{n}{j}(-1)^{n-j}2^j b_j . \qquad (2.6)$$

**Example 6**. Here is one very simple demonstration how the corollary works. Let $a_n = (-1)^{n-1}$. Then we have for $n \geq 0$

$$\sum_{k=0}^{n}\binom{n}{k}(-1)^{k-1} = \begin{cases} -1, n=0 \\ 0, n>0 \end{cases}$$

and from (1.8) we find $\nabla^m b_n = 0$ for $m \neq n$ and $\nabla^n b_n = (-1)^{n-1}$. Therefore, from (2.5),

$$\sum_{k=0}^{n}\binom{n}{k}(-1)^{k-1} x^k = (-1)^{n-1}(x-1)^n = -(1-x)^n . \qquad (2.7)$$

Of course, this identity follows immediately from the binomial formula.

## 3. Identities with special numbers

This section contains some new identities for products of harmonic, Bernoulli and Fibonacci numbers. We start with the following lemma:

**Lemma 3.** For any two integers $1 \leq m \leq n$,

$$\sum_{k=m}^{n} \binom{n}{k}\binom{k}{m}\frac{(-1)^k}{k} = \frac{(-1)^m}{m} . \tag{3.1}$$

The proof is given in Section 5.

**Example 7.** Let $a_k = \frac{(-1)^k}{k}$, $k \geq 1$. Them for any $1 \leq m \leq n$ we have from the above lemma and from Lemma 1,

$$\sum_{k=1}^{n} \binom{n}{k}\binom{k}{m} a_k = \binom{n}{m} \nabla^m b_n = \frac{(-1)^m}{m} . \tag{3.2}$$

Then from Theorem 1, the symmetric version, equation (1.5)

$$\sum_{k=1}^{n} \binom{n}{k}\frac{(-1)^k c_k}{k} = \sum_{m=1}^{n} \frac{d_m}{m} . \tag{3.3}$$

This property of the binomial transform was discussed in [1]. It is true for any two sequences $\{c_k\}, \{d_k\}, k \geq 1$ related by

$$d_n = \sum_{k=0}^{n} \binom{n}{k}(-1)^k c_k, \quad c_n = \sum_{k=0}^{n} \binom{n}{k}(-1)^k d_k ,$$

and the factor $(-1)^k$ can be replaced by $(-1)^{k-1}$. For the transform (1.1) the property is

$$\sum_{k=1}^{n} \binom{n}{k}\frac{a_k}{k} = \sum_{m=1}^{n} \frac{b_m}{m} .$$

**Example 8**. This example is related to the previous one. Let

$$H_n = 1 + \frac{1}{2} + ... + \frac{1}{n}, \quad H_0 = 0, \quad n \geq 0,$$

be the harmonic numbers. Then we have (see [4])

$$\sum_{k=1}^{n} \binom{n}{k} \frac{(-1)^{k-1}}{k} = H_n . \tag{3.4}$$

According to Lemma 1,

$$\sum_{k=m}^{n} \binom{n}{k}\binom{k}{m} \frac{(-1)^{k-1}}{k} = \binom{n}{m} \nabla^m H_n \tag{3.5}$$

and then from Lemma 3, for $1 \leq m \leq n$,

$$\binom{n}{m} \nabla^m H_n = \frac{(-1)^{m-1}}{m} \tag{3.6}$$

and when $m=0$ we have $\binom{n}{0} \nabla^0 H_n = H_n$.

Therefore, from (3.4) and Corollary 2, by separating the first term in the second sum,

$$\sum_{k=1}^{n} \binom{n}{k} \frac{(-1)^{k-1} x^k}{k} = \sum_{m=0}^{n} \binom{n}{m} \nabla^m H_n (x-1)^m = H_n + \sum_{k=1}^{n} \frac{(-1)^{k-1}(x-1)^k}{k} ,$$

that is,

$$\sum_{k=1}^{n} \binom{n}{k} \frac{(-1)^{k-1} x^k}{k} = H_n - \sum_{k=1}^{n} \frac{(1-x)^k}{k} \tag{3.7}$$

This equation also follows from (2.7) and (3.3).

Formula (2.5) can be used to evaluate the iterated differences $\nabla^m b_n$ when the LHS in (2.5) is known

**Example 9**. By inversion in (3.4) we have

$$\sum_{k=0}^{n} \binom{n}{k} (-1)^{k-1} H_k = \frac{1}{n} . \tag{3.8}$$

The version with "$x$" was presented in [2], that is,

$$\sum_{k=0}^{n} \binom{n}{k} (-1)^{k-1} x^k H_k = \frac{1}{n} + \frac{1-x}{n-1} + \frac{(1-x)^2}{n-2} + ... + \frac{(1-x)^{n-2}}{2} + \frac{(1-x)^{n-1}}{1} - (1-x)^n H_n .$$

Comparing this to (2.5) we conclude that

$$\binom{n}{m} \nabla^m \frac{1}{n} = \frac{(-1)^m}{n-m} \quad \text{when } 0 \leq m < n, \text{ and } \nabla^n \frac{1}{n} = (-1)^{n-1} H_n. \tag{3.9}$$

From Example 9 and Theorem 1 we derive the following most interesting result:

**Corollary 3**. Let $\{c_k\}$ and $\{d_k\}$ be any two sequences related as in (1.2). Then

$$\sum_{k=0}^{n} \binom{n}{k} (-1)^{k-1} H_k c_k = (-1)^{n-1} H_n d_n + \sum_{m=0}^{n-1} \frac{(-1)^m d_m}{n-m}. \tag{3.10}$$

For the proof we take $a_k = (-1)^{k-1} H_k$ and then in view of (1.1) and (3.8) we have $b_n = \frac{1}{n}$. The rest follows from (3.9) and Theorem 1.

To show Corollary 3 in action we shall give several examples.

**Example 10.** Applying property (3.3) to equation (3.8) we find

$$\sum_{k=1}^{n} \binom{n}{k} (-1)^{k-1} \frac{H_k}{k} = 1 + \frac{1}{2^2} + \ldots + \frac{1}{n^2}$$

With the notation

$$H_n^{(2)} = 1 + \frac{1}{2^2} + \ldots + \frac{1}{n^2}, \quad H_0^{(2)} = 0,$$

we obtain by inversion $(n \geq 1)$

$$\frac{H_n}{n} = \sum_{k=0}^{n} \binom{n}{k} (-1)^{k-1} H_k^{(2)},$$

and (3.10) yields (with $c_k = \frac{H_k}{k}$, $d_m = (-1)^{m-1} H_m^{(2)}$, related as in (1.2)),

$$\sum_{k=0}^{n} \binom{n}{k} (-1)^{k-1} \frac{H_k^2}{k} = H_n H_n^{(2)} - \sum_{m=0}^{n-1} \frac{H_m^{(2)}}{n-m}. \tag{3.11}$$

**Example 11**. In the same way, starting from the identity (see [1], equation (16))

$$\frac{H_n}{n+1} = \sum_{k=0}^{n} \binom{n}{k} (-1)^{k-1} \frac{H_k}{k+1},$$

and taking $c_k = \dfrac{H_k}{k+1}$, $d_m = \dfrac{(-1)^{m-1} H_m}{m+1}$ in (3.10) we obtain

$$\sum_{k=0}^{n} \binom{n}{k} (-1)^{k-1} \frac{H_k^2}{k+1} = \frac{H_n^2}{n+1} - \sum_{m=0}^{n-1} \frac{H_m}{(n-m)(m+1)}. \tag{3.12}$$

Now we show identities involving other special numbers.

**Example 12.** By inversion in (2.3),

$$n^\alpha = \sum_{k=0}^{n} \binom{n}{k} k! S(\alpha, k)$$

and therefore, from (3.10) with $c_k = k^\alpha$ and $d_m = m! S(\alpha, m)$,

$$\sum_{k=0}^{n} \binom{n}{k} (-1)^{k-1} H_k k^\alpha = (-1)^{n-1} n! H_n S(\alpha, n) + \sum_{m=0}^{n-1} \frac{(-1)^m m! S(\alpha, m)}{n-m}. \tag{3.13}$$

**Example 13.** For the Fibonacci numbers $F_n$ it is known that the following two binomial identities are true

$$F_n = \sum_{k=0}^{n} \binom{n}{k} (-1)^{k-1} F_k, \tag{3.14}$$

$$F_{2n} = \sum_{k=0}^{n} \binom{n}{k} F_k.$$

From here and (3.10) we derive two new identities involving products of harmonic and Fibonacci numbers

$$\sum_{k=0}^{n} \binom{n}{k} (-1)^{k-1} H_k F_k = H_n F_n - \sum_{m=0}^{n-1} \frac{F_m}{n-m}, \tag{3.15}$$

by taking $c_k = F_k$, $d_m = (-1)^{m-1} F_m$. Also

$$\sum_{k=0}^{n} \binom{n}{k} (-1)^{k-1} H_k F_{2k} = (-1)^{n-1} H_n F_n + \sum_{m=0}^{n-1} \frac{(-1)^m F_m}{n-m} \tag{3.16}$$

with $c_k = F_{2k}$, $d_m = F_m$.

**Example 14**. Here we use the Bernoulli numbers $B_n$ defined by the generating function

$$\frac{t}{e^t - 1} = \sum_{n=0}^{\infty} B_n \frac{t^n}{n!}, \quad |t| < 2\pi.$$

For the Bernoulli numbers it is known that

$$(-1)^n B_n = \sum_{k=0}^{n} \binom{n}{k} B_k. \qquad (3.17)$$

From (3.10) with $c_k = (-1)^k B_k$, $d_m = B_m$ we obtain the identity

$$\sum_{k=0}^{n} \binom{n}{k} H_k B_n = (-1)^n H_n B_n - \sum_{m=0}^{n-1} \frac{(-1)^m B_m}{n-m}. \qquad (3.18)$$

**Example 15**. In this last example related to (3.10) we use the Laguerre polynomials

$$L_n(x) = \frac{e^x}{n!} \left(\frac{d}{dx}\right) (x^n e^{-x})$$

which satisfy the identity

$$L_n(x) = \sum_{k=0}^{n} \binom{n}{k} \frac{(-x)^k}{k!}.$$

Here (3.10) provides the curious formula ( $c_k = L_k(x)$ and $d_m = \frac{(-x)^m}{m!}$ )

$$\sum_{k=0}^{n} \binom{n}{k} (-1)^k H_k L_k(x) = \frac{x^n}{n!} H_n - \sum_{m=0}^{n-1} \frac{x^m}{m!(n-m)}. \qquad (3.19)$$

Next we turn again to the sequence of Fibonacci numbers defined by the recurrence $F_n = F_{n-1} + F_{n-2}$, and starting with $F_0 = 0$, $F_1 = 1$. We can extend the sequence $F_n$ for negative indices by using the equation $F_{n-2} = F_n - F_{n-1}$. Thus we come to the negatively indexed Fibonacci numbers, where $F_{-n} = (-1)^{n+1} F_n$, $n \geq 0$. Computing the backward differences we find

$$\nabla F_n = F_n - F_{n-1} = F_{n-2},$$

$$\nabla^2 F_n = F_{n-2} - F_{n-3} = F_{n-4},$$

etc. Obviously, $\nabla^m F_n = F_{n-2m}$ and this is true for any non-negative integer $m$. Now we can formulate the desired result:

**Corollary 4.** For any pair of sequences $\{c_k\}$ and $\{d_k\}$ as in (1.2), and every non-negative integer $n$ we have

$$\sum_{k=0}^{n}\binom{n}{k}(-1)^{k-1}F_k\, c_k = \sum_{m=0}^{n}\binom{n}{m}d_m F_{n-2m} \ . \tag{3.20}$$

For the proof we use (1.4) in Theorem 1 with $a_k = (-1)^{k-1}F_k$ and $b_n = F_n$ (see (3.14)).

Formula (3.20) can be used in the same way as (3.10) to generate various new identifies by choosing different sequences $\{c_k\}$. For illustration we provide the following example:

**Example 16.** Choosing $c_n = (-1)^n B_n$ and $d_k = B_k$, where $B_n$ are the Bernoulli numbers, we obtain from (3.17) and (3.20) an identity connecting Bernoulli and Fibonacci numbers

$$\sum_{k=0}^{n}\binom{n}{k}B_k F_k = -\sum_{m=0}^{n}\binom{n}{m}B_m F_{n-2m} \tag{3.21}$$

or,

$$\sum_{k=0}^{n}\binom{n}{k}B_k (F_k + F_{n-2k}) = 0 \ .$$

The Lucas numbers $L_n$ satisfy the same recurrence $L_n = L_{n-1} + L_{n-2}$ as the Fibonacci numbers and for them binomial identities like (3.14) hold too. Therefore, a property similar to (3.20) is also true for the Lucas numbers.

## 4. Some variations

**Remark 2.** If the binomial transform is defined by the formula

$$\sum_{k=0}^{n}\binom{n}{k}(-1)^k a_k = b_n \ ,$$

then (2.5) takes the form

$$\sum_{k=0}^{n}\binom{n}{k}a_k(-x)^k = \sum_{m=0}^{n}\binom{n}{m}\nabla^m b_n (x-1)^m \tag{4.1}$$

or

$$\sum_{k=0}^{n}\binom{n}{k}a_k x^k = \sum_{m=0}^{n}\binom{n}{m}(-1)^m \nabla^m b_n (x+1)^m.$$

**Remark 3**. Another expression for the coefficients

$$C(n,m) = \binom{n}{m}\nabla^m b_n$$

can be written in terms of the Stirling numbers of the first kind $s(m, j)$. A good reference for these numbers is the book [6].

Suppose the sequence $\{b_n\}$ is the binomial transform of the sequence $\{a_n\}$ as in (1.1). Then

$$C(n,m) = \frac{1}{m!}\sum_{j=0}^{m} s(m, j)(n\nabla)^j b_n \qquad (4.2)$$

*Proof.* We have the representation

$$\binom{k}{m} = \frac{1}{m!}\sum_{j=0}^{m} s(m, j) k^j$$

(see [6]) and from here and (1.7)

$$C(n,m) = \sum_{k=0}^{n}\binom{n}{k}\binom{k}{m}a_k = \frac{1}{m!}\sum_{k=0}^{n}\binom{n}{k}\left\{\sum_{j=0}^{m} s(m, j) k^j a_k\right\}$$

$$= \frac{1}{m!}\sum_{j=0}^{m} s(m, j)\left\{\sum_{k=0}^{n}\binom{n}{k} k^j a_k\right\}$$

Therefore, in view of (1.3) we come to (4.2).

## 5. Proofs

Here we prove the three lemmas and Theorem 1.

*Proof of Lemma 1.* When $n=0$ this is obviously true. Take any integers $n \geq 1$. We shall do induction on $1 \leq m \leq n$. Suppose the identity is true for some $m < n$. We shall prove it for $m+1$.

The LHS then becomes (using (1.3) with $p=1$ in the second equality)

$$\sum_{k=0}^{n}\binom{n}{k}k(k-1)...(k-m+1)(k-m)a_k = \sum_{k=0}^{n}\binom{n}{k}k\{k(k-1)...(k-m+1)a_k\}$$

$$-m\sum_{k=0}^{n}\binom{n}{k}k(k-1)...(k-m+1)a_k = n\nabla\left\{m!\binom{n}{m}\nabla^m b_n\right\} - m\left\{m!\binom{n}{m}\nabla^m b_n\right\}$$

$$= m!\left\{n\binom{n}{m}\nabla^m b_n - n\binom{n-1}{m}\nabla^m b_{n-1} - m\binom{n}{m}\nabla^m b_n\right\}$$

$$= m!\left\{(n-m)\binom{n}{m}\nabla^m b_n - n\binom{n-1}{m}\nabla^m b_{n-1}\right\}$$

$$= m!\left\{\frac{n!}{(n-m-1)!m!}\nabla^m b_n - \frac{n!}{(n-m-1)!m!}\nabla^m b_{n-1}\right\}$$

$$= m!\left\{\frac{n!}{(n-m-1)!m!}\nabla^{m+1} b_n\right\} = (m+1)!\binom{n}{m+1}\nabla^{m+1} b_n .$$

*Proof of Lemma 2.* From Lemma 1 and the inversion formula for the binomial transform

$$\binom{n}{m}\nabla^m b_n = \frac{1}{m!}\sum_{k=0}^{n}\binom{n}{k}k(k-1)...(k-m+1)a_k = \sum_{k=0}^{n}\binom{n}{k}\binom{k}{m}a_k$$

$$= \sum_{k=0}^{n}\binom{n}{k}\binom{k}{m}\left\{\sum_{j=0}^{k}\binom{k}{j}(-1)^{k-j}b_j\right\} = \sum_{j=0}^{n}(-1)^j b_j\left\{\sum_{k=0}^{n}\binom{n}{k}\binom{k}{m}\binom{k}{j}(-1)^k\right\}$$

$$= \sum_{j=0}^{n}\binom{n}{j}\binom{j}{n-m}(-1)^{n-j}b_j$$

by using the identity ([7], p.15)

$$\sum_{k=0}^{n}\binom{n}{k}\binom{k}{m}\binom{k}{j}(-1)^k = (-1)^n\binom{n}{j}\binom{j}{n-m} .$$

*Proof of Theorem1.*

$$\sum_{k=0}^{n}\binom{n}{k}a_k c_k = \sum_{k=0}^{n}\binom{n}{k}a_k\left\{\sum_{m=0}^{k}\binom{k}{m}d_m\right\}$$

$$= \sum_{m=0}^{n} d_m \left\{ \sum_{k=0}^{n} \binom{n}{k}\binom{k}{m} a_k \right\} = \sum_{m=0}^{n} \binom{n}{m} d_m \nabla^m b_n$$

according to Lemma 1.

*Proof of Lemma 3.*

The starting point of this proof is the identity (2.2) with $1 \leq j \leq n$. We divide both sides by $x$ and integrate from 0 to 1. This yields

$$\sum_{k=j}^{n} \binom{n}{k}\binom{k}{j} \frac{(-1)^k}{k} = (-1)^j \binom{n}{j} \int_0^1 x^{j-1}(1-x)^{n-j} dx$$

$$= (-1)^j \binom{n}{j} B(n-j+1, j) = \frac{(-1)^j}{j} \ .$$

The evaluation of the integral is from table [5], namely, this is entry 3.191 (3). Here $B(x, y)$ is Euler's Beta function,

$$B(x, y) = \frac{\Gamma(x)\Gamma(y)}{\Gamma(x+y)} \ .$$